\documentclass[12pt]{article}
\usepackage{epsfig,makeidx,amsmath,amsfonts,latexsym,subfigure}
\usepackage{times,pifont}
\usepackage{epsf}
\usepackage{graphics}
\newtheorem{theorem}{Theorem}[section]
\newtheorem{proposition}[theorem]{Proposition}
\newtheorem{lemma}[theorem]{Lemma}
\newtheorem{corollary}[theorem]{Corollary}
\newtheorem{definition}[theorem]{Definition}

\newenvironment{proof}{{\em Proof. }}{\hfill $\Box$ \vspace{1em}}

\newcommand{\E}{{\mathbb E}}

\newcommand{\Var}{{\mathbb V}{\rm ar}}
\newcommand{\Cov}{{\mathbb C}{\rm ov}}
\newcommand{\RR}{{\mathbb R}}
\newcommand{\CC}{{\mathbb C}}

\newcommand{\SG}{{\mathbb S}}

\newcommand{\Pro}{{\mathbb P}}

\newcommand{\taum}{\tau_{{\rm mix}}}

\newcommand{\n}{\|}
\newcommand{\bi}{\begin{itemize}}
\newcommand{\ei}{\end{itemize}}
\newcommand{\be}{\begin{enumerate}}
\newcommand{\ee}{\end{enumerate}}
\newcommand{\ra}{\rightarrow}

\newcommand{\ep}{\epsilon}

\newcommand{\iy}{\infty}

\newcommand{\beq}{\begin{equation}}
\newcommand{\eeq}{\end{equation}}
\newcommand{\beqa}{\begin{eqnarray*}}
\newcommand{\eeqa}{\end{eqnarray*}}
\newcommand{\btm}{\begin{theorem}}
\newcommand{\etm}{\end{theorem}}
\newcommand{\bpf}{\begin{proof}}
\newcommand{\epf}{\end{proof}}
\newcommand{\bla}{\begin{lemma}}
\newcommand{\ela}{\end{lemma}}
\newcommand{\bdn}{\begin{definition}}
\newcommand{\edn}{\end{definition}}
\newcommand{\bpn}{\begin{proposition}}
\newcommand{\epn}{\end{proposition}}
\newcommand{\bcy}{\begin{corollary}}
\newcommand{\ecy}{\end{corollary}}

\begin{document}

\title{Card-cyclic-to-random shuffling with relabeling}
\author{Johan Jonasson\thanks{Chalmers University of Technology and University of Gothenburg} \thanks{Research supported by the Knut and Alice Wallenberg Foundation}}
\maketitle

\begin{abstract}
The card-cyclic-to-random shuffle is the card shuffle where the $n$ cards are labeled $1,\ldots,n$ according to their starting positions. Then the cards are mixed by first picking card $1$ from the deck and reinserting it at a uniformly random position, then repeating for card $2$, then for card $3$ and so on until all cards have been reinserted in this way. Then the procedure starts over again, by first picking the card with label $1$ and reinserting, and so on.
Morris, Ning and Peres \cite{MNP} recently showed that the order of the number of shuffles needed to mix the deck in this way is $n\log n$.
In the present paper, we consider a variant of this shuffle with relabeling, i.e.\ a shuffle that differs from the above in that after one round, i.e.\ after all cards have been reinserted once, we relabel the cards according to the positions in the deck that they now have. The relabeling is then repeated after each round of shuffling.
It is shown that even in this case, the correct order of mixing is $n\log n$.

\end{abstract}

\noindent{\em Short title :\/ Card-cyclic-to-random shuffling } 

\section{Introduction}

The subject of mixing times for Markov chains an important and exciting research field that has attracted a lot of attention in recent decades.
An outstanding subclass of Markov chains that has been studied extensively is card shuffling, i.e.\ Markov chains on the symmetric group $\SG_n$ of permutations of
$n$ items that one can think of as the cards of a deck.

One of the early card shuffles to be studied was the random transpositions shuffle, where each step of the shuffle is made by picking two cards uniformly and independently at random and then swapping them.
It was shown by Diaconis and Shahshahani \cite{DSh} that the mixing time of this shuffle has a sharp threshold at $\frac12 n\log n$ shuffles.
It is easy to see that at least order of $n \log n$ shuffles is required, since, by the coupon collector's problem, it takes this order of shuffles until most cards have been touched at all.
Closely related to the random transpositions shuffle is the top-to-random shuffle where at each step the card presently in position one is moved to a uniform random position.
The sharp threshold for this shuffle is $n\log n$ and again it is easy to see that at least order of $n\log n$ steps is required for mixing, for similar reasons.

In recent years some more systematic variants of these shuffles have been proposed and analyzed. Mossel, Peres and Sinclair \cite{MPS} and Saloff-Coste and Zuniga \cite{SZ} analyzed the
cyclic-to-random shuffle, where at time $t$ the card presently in position $t$ mod $n$ is swapped with a uniformly random card.
Clearly at least once per $n$ steps, each card will be touched and one of the interesting questions about this shuffle was if $O(n)$ shuffles is also sufficient to mix the whole deck.
The answer turns out to be negative; indeed the mixing time is still of order $n\log n$.
Pinsky \cite{Pinsky} later introduced the card-cyclic-to-random transpositions shuffle (CCR shuffle), where at time $t$ the card with {\em label} $t$ mod $n$ (i.e.\ the card that started out in position $t$ mod $n$) is moved to a uniformly random position.
Again it is obvious that every card will be touched once every $n$ steps and again one main question was if this way of systematically randomizing the cards, suffices to mix the whole deck in $O(n)$, or at least $o(n\log n)$, steps.
Again the answer turns out to be negative; Morris, Ning and Peres \cite{MNP} prove that $n\log n$ is still the correct order.
In this paper we investigate the {\em card-cyclic-to-random shuffle with relabeling} (the CCRR shuffle for short). For $k=1,2,\ldots$ let {\em round} $k$ consist of steps
$kn+1,kn+2,\ldots,kn+n$ of shuffling. The CCRR shuffle is the shuffle that is exactly as the card-cyclic-to-random shuffle for the first round.
After that however, the cards are relabeled $1,\ldots,n$ according to their positions after the first round. Next a new round of CCR shuffling is carried out according to the new labels.
After that the cards are relabeled again and a new round of CCR is done, and so on.
The main result of this paper is that relabeling does not help to speed up mixing either, at least not more than by a constant.

\btm \label{ta}
The mixing time of the card-cyclic-to-random transpositions with relabeling is of order $n\log n$.
\etm

Here, the mixing time is given by
\[\taum := \min\{t: \n \Pro(X_t \in \cdot) - \pi\n_{TV} \leq \frac14\}\]
where $X_t \in \SG_n$ is the state of the deck of cards after $t$ steps of shuffling, $\pi$ is the uniform distribution on $\SG_n$ and $\n \cdot\n_{TV}$ is the total variation norm, given
in general by
\[\n\mu\n_{TV} := \frac12\sum_{x \in S}|\mu(x)| = \max\{\mu(A): A \subset S\}\]
for a signed measure $\mu$ on a finite space $S$.

\section{Proof of the main result}
For the upper bound on $\taum$, it suffices to note that the proof in \cite{MNP} for the CCR shuffle goes through exactly as it stands there. Hence we will focus entirely on the lower bound.
The idea of the proof of the lower bound draws on the idea behind Wilson's technique introduced in \cite{Wilson1} and \cite{Wilson2}, namely to use an eigenvector of the transition matrix for the movement of a single card to build a test function. However since estimating the variance of the test function will in fact be quite simple here, we will not need
Wilson's Lemma explicitly.
A rough outline of the proof is
\be
\item Show that the position of a given card after one round of CCRR is determined, up to a random term of order $\sqrt{n}$, by where it was reinserted.
\item In the light of 1, study the idealized motion of a single card which is a deterministic function of where it was reinserted.
\item Show that the transition matrix for one round of idealized single card motion has a spectral gap bounded away from $1$.
\item Use the eigenvector corresponding to the second eigenvalue to construct a test statistic.
\item Estimate, using 1, the expectation and variance of the test statistic applied to the CCRR shuffle and establish the lower bound using Chebyshev's inequality.
\ee

Because of the cyclic structure of the shuffle, the movement of a single card is not time-homogenous if we consider individual steps of the shuffle. However in terms of {\em rounds}, the movement of a given card is indeed a time-homogenous Markov chain.
Let $A=A(n)$ denote the transition matrix of this chain on $n$ cards.
It turns out that when analyzing this Markov chain, it is convenient to denote the possible positions a card can have in the deck as $1/n,2/n,3/n,\ldots,1$ (instead of the usual $1,2,\ldots,n$).
Write $Q_n := [n]/n = \{1/n,2/n,\ldots,1\}$ for the set of positions.
As usual, we will identify a card with its starting position, i.e.\ when we speak of card $a$, $a \in Q_n$, we are considering the card that starts in position $a$.
Since this is the $na$'th card from the top in the starting order of the deck, we may sometimes also speak of this card as card $na$.

It is difficult to come up with a closed-form expression for $A$, but the action of $A$ can be probabilistically described as follows. Consider a card that starts a round in position $a \in Q_n$.
Let us refer to the cards $1/n,\ldots,a-1/n$ as white cards and to the cards $a+1/n,\ldots,1$ as black cards.
Now in a first stage the $na-1$ white cards are sequentially picked out and reinserted at independent uniform positions. During this stage a certain number of cards will be reinserted above card $a$ in the deck whereas the others will be uniformly spread out among the black cards below card $a$. The cards that in this stage end up above card $a$ will form a well-mixed layer of white cards. Note that during stage 1, card $a$ will move gradually higher up in the deck.
(Here we say that if $a<b$, then position $a$ is higher up than, or above, position $b$.)

Next, after stage 1, card $a$ itself is picked out and reinserted at a uniformly random position $U = U_n \in Q_n$; this is stage 2.
In the third and final stage, the black cards are picked out and reinserted. If card $a$ was reinserted in the white layer at the top, then card $a$ will move gradually down the deck during the whole of this stage, whereas if not, then stage 3 divides into the two sub-stages where in the first of these, stage 3a, the black cards above card $a$ are reinserted and $a$ moves upwards and in the second, stage 3b, the black cards below card $a$ are reinserted and $a$ moves down the deck.

Even though we will not need the exact distribution of where card $a$ ends up under this procedure, we will still need a good approximate control.
The following two lemmas will be useful for that.

\bla \label{la}
Let the sequence $Y_0,Y_1,\ldots,Y_{n(1-a)}$ be recursively defined by $Y_0=a \in Q_n$ and $Y_{t+1}=Y_t+1/n$ with probability $Y_t$ and $Y_{t+1}=Y_t$ with probability $1-Y_t$ (where these events are conditionally independent of $Y_0,Y_1,\ldots,Y_{t-1}$ given $Y_t$).
Then
\[\E[Y_t] = \left(1+\frac{1}{n}\right)^t a\]
and
\[\Var(Y_{t+1}) \leq \frac{a}{n^2}\sum_{j=t}^{2t}\left(1+\frac{1}{n}\right)^j - \frac{a^2}{n^2}t\left(1+\frac{1}{n}\right)^{2t}.\]
In particular,
for all $t$,
\[\Var(Y_t) < \frac25 n^{-1}.\]
\ela

\bpf
By conditioning on $Y_{t}$ we get that
\[\E[Y_{t+1}] = \E\left[Y_t\left(Y_t+\frac1n\right)+\left(1-Y_t\right)Y_t\right] = \E\left[\left(1+\frac{1}{n}\right)Y_t\right]\]
which proves the expression for the expectation.
For the variance part, write $v_t:=\Var(Y_t)$. Then $v_0=0$ and recursively
\[\Var(Y_{t+1}) = \E[\Var(Y_{t+1}|Y_t)] + \Var(\E[Y_{t+1}|Y_t]).\]
By definition of the $Y_t$'s, $\Var(Y_{t+1}|Y_t) = Y_t(1-Y_t)/n^2$ and by the above $\E[Y_{t+1}|Y_t] = (1+1/n)Y_t$.
For for first term we have
\beqa
\E[\Var(Y_{t+1}|Y_t)] &=& \frac{\E[Y_t]}{n^2}-\frac{\E[Y_t^2]}{n^2} \\
&=& \frac{(1+1/n)^t y_0}{n^2}-\frac{1}{n^2}\left(v_t+a^2\left(1+\frac{1}{n}\right)^{2t}\right).
\eeqa
Adding the second term and writing $c:=1+1/n$ gives
\begin{eqnarray}
v_{t+1} &=& \left(c^2-\frac{1}{n^2}\right)v_t+\frac{c^t a}{n^2}\left(1-c^t a\right) \label{eeu} \\
&<& c^2v_t+\frac{1}{n^2}c^t a - \frac{1}{n^2}c^{2t}a^2 \nonumber.
\end{eqnarray}
This recursion is readily solved and gives
\[v_{t+1} < \frac{a}{n^2}\sum_{j=t}^{2t}c^j - \frac{a^2}{n^2}tc^{2t}.\]
By (\ref{eeu}), $v_t$ is increasing in $t$, so we get an upper bound on plugging in $t=n(1-a)$ on the right hand side and then get
\beqa
v_{t+1} &<& \frac{ea}{n}\left(e^{1-2a}-e^{-a}-a(1-a)e^{1-2a}\right) \\
&<& \frac25 n^{-1},
\eeqa
where the second inequality from standard optimization over $a$.

\epf

\bla \label{lcontractive}
Let $X \in L^2(\RR)$ be a random variable and $f:\RR \ra \RR$ be contractive, i.e.\
$|f(x)-f(y)| \leq |x-y|$ for all $x,y \in \RR$.
Then
\[\Var(f(X)) \leq \Var(X).\]
\ela

\bpf
Let $X_1$ and $X_2$ be two independent copies of $X$. Then
\beqa
\Var(X) &=& \frac12\Var(X_1-X_2)\\
&=& \frac12 \E[|X_1-X_2|^2] \\
&\geq& \frac12 \E[|f(X_1)-f(X_2)|^2]\\
&=& \Var(f(X)).
\eeqa
\epf

\begin{figure}
\begin{center}
\includegraphics[trim = 35mm 95mm 40mm 100mm, clip, width=0.7\textwidth]{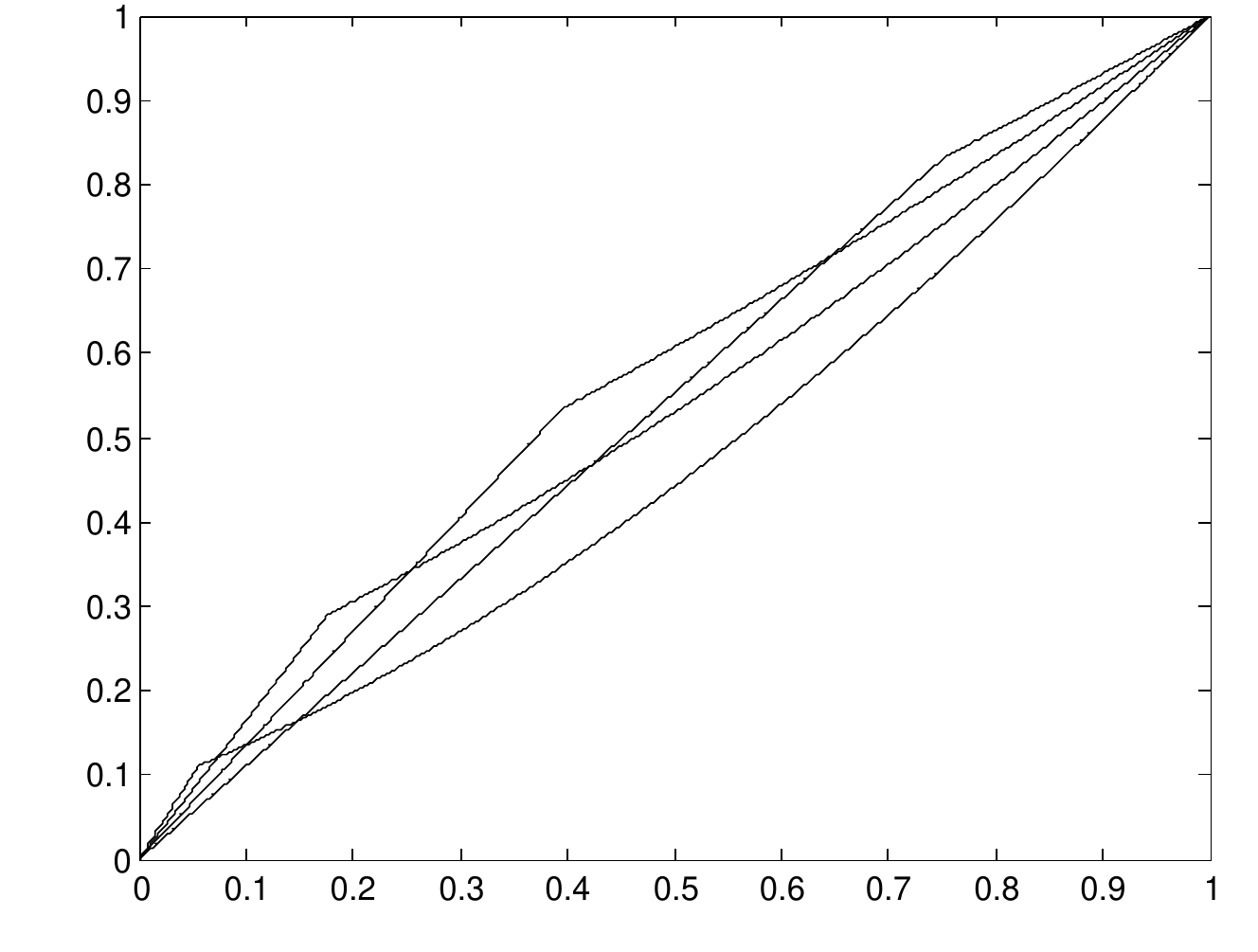}
\caption{The function $G_b(u)$ for $b=0.3,0.5,0.7,0.9$. The smaller the $b$, the larger the ascent at the origin.} \label{fa}
\end{center}
\end{figure}

Let $Z=Z_n$ be the position that card $a=a_n\in Q_n$ ends up in after one round of shuffling.
For each $b \in [0,1]$, define $G_b:[0,1] \ra [0,1]$ as
\beq \label{eideal}
G_b(u) = \left\{ \begin{array}{ll} e^{1-b}u, & u \leq u_0(b) := 1-(1-b)e^b \\ e^{e^{-b}(1-u)}-(1-u)e^{1-b}, & u>u_0(b) \end{array} \right.
\eeq
Note that $G_b(u)$ is continuous in $(b,u)$ and for each $b$, $G_b$ is differentiable for $u \neq u_0(b)$.
See Figure \ref{fa} to see a plot of $G_b$ for a few different $b$.
The functions $G_b$ play a central r\^ole in the following lemma, which gives control over the asymptotic distribution, expectation and variance of $Z$ given $U$.
The limiting distribution is known and due to Pinsky, see Theorem 4 of \cite{Pinsky}. Since we will also need a quantified bound on the variance, we will for self containedness, reprove the result below.

\bla \label{lappr}
For all $u \in Q_n$,
\beq \label{ei}
\E[Z|U=u] \in \left(1 \pm \frac{2}{n}\right)G_a(u).
\eeq
and
\[\Var(Z|U=u) < 9n^{-1},\]
or equivalently, writing $Z=\E[Z|U]+D$,
\beq \label{ea}
\Var(D) < 9n^{-1}.
\eeq
In particular, if $a_n \ra a \in [0,1]$, then
\[Z_n \stackrel{d}{\ra} G_a(V)\]
where $V$ is uniform on $[0,1]$.
\ela

\bpf
Let $W = W_n \in Q_n$ be such that $nW$ is the number of white cards that go to the top layer of white cards in stage 1. We will start by estimating the variance of $Z$ given $W=w$ and $U=u$.
If $u \leq w$, so that stage 2 moves card $a$ to the top white layer, then by Lemma \ref{la}
\[\E[Z|U=u,W=w] = \left(1+\frac{1}{n}\right)^{n(1-a)}u\]
and
\[\Var(Z|U=u,W=w)<\frac25 n^{-1}.\]
The case $u > w$ takes some more work. In order to not overly burden the notation, we will until further notice, point out the conditioning on $U=u$ and/or $W=w$ by writing
indexes $u$ and/or $w$ at the conditional expectations and variances.

Let $S \in Q_n$ be such that $nS$ is the number of black cards in positions $w+1/n,\ldots,u-1/n$; these are the black cards that get reinserted in stage 3a.
Let $nN$ be the number of cards below card $a$ after these $nS$ black cards have been reinserted. Note that when the first part of stage 3 starts, then card $a$ is in $u$ and
at that point, the number of black cards below $a$ is $n(1-a-S)$; these are the ones that will get reinserted in the stage 3b.
Hence by Lemma \ref{la},
\[\E_{u,w}[Z|N,S] = \left(1+\frac{1}{n}\right)^{n(1-a-S)}(1-N)\]
and
\[\E_{u,w}[N|S] = \left(1+\frac{1}{n}\right)^{nS}(1-u).\]
Hence
\beqa
\E_{u,w}[Z|S] &=& (1-\E_{u,w}[N|S])\left(1+\frac{1}{n}\right)^{n(1-a-S)} \\
&=& \left(1+\frac{1}{n}\right)^{n(1-a-S)} - \left(1+\frac{1}{n}\right)^{n(1-a)}(1-u).
\eeqa
It follows that
\beq \label{eh}
\E_{u,w}[Z] = \E_{u,w}\left[ \left(1+\frac{1}{n}\right)^{n(1-a-S)} \right] - \left(1+\frac{1}{n}\right)^{n(1-a)}(1-u).
\eeq
We also get that
\beqa
\Var_{u,w}(Z|S) &=& \E_{u,w}\left[\Var_{u,w}(Z|N,S)|S\right] + \Var_{u,w}\left(\E_{u,w}[Z|N,S]|S\right) \\
&<& \frac25 n^{-1} + \Var_{u,w} \left( \left(1+\frac{1}{n}\right)^{n(1-a-S)}(1-N)\Big|S\right) \\
&=& \frac25 n^{-1} + \left(1+\frac{1}{n}\right)^{2n(1-a-S)}\frac25 n^{-1} \\
&<& \frac25 (1+e^2)n^{-1}.
\eeqa
Therefore
\beqa
\Var_{u,w}(Z) &=& \E_{u,w}\left[\Var_{u,w}(Z|S)\right] + \Var_{u,w}\left(\E_{u,w}[Z|S]\right) \\
&<& \frac25 (1+e^2)n^{-1}+\Var_{u,w}\left( \left(1+\frac{1}{n}\right)^{n(1-a-S)}\right) \\
&\leq& \frac25 (1+e^2)n^{-1}+\Var_{u,w}(eS) \\
&<& \left(\frac25+\frac{13}{20}e^2\right)n^{-1}
\eeqa
by Lemma \ref{lcontractive}, since the map $S \ra (1/e)(1+1/n)^{n(1-a-S)}$ is contractive and $nS$, given $U=u$ and $W=w$, is hypergeometric with variance at most $n/4$.

Now bring back the conditioning on $W$ into ordinary notation.
What we have just shown is among other things, that $\Var_u(Z|W) < Cn$ with $C:=2/5+13e^2/20$.
Thus $\Var_u(Z) = \E_u[\Var_u(Z|W)]+\Var_u(\E_u[Z|W]) < Cn+\Var_u(\E_u[Z|W])$.
However, by (\ref{eh})
\[\left| \E_u\left[Z|W=w\right]-E_u\left[Z|W=w-1/n\right] \right|\]
\[= \left| \E_u\left[\left(1+\frac{1}{n}\right)^{n(1-a-S)}\Big|W=w\right]- \E_u\left[\left(1+\frac{1}{n}\right)^{n(1-a-S)}\Big|W=w-1/n\right] \right|\]
\[\leq \E_u \left[ \left( 1+\frac{1}{n} \right)^{n(1-a-S)} \left( \left( 1+\frac{1}{n} \right)-1 \right) \Big| W=w \right] \leq \frac{e}{n},\]
where the first inequality uses that the conditional distributions of $nS$ given $W=w$ and $W=w-1/n$ respectively, can easily be coupled so that the realizations do not differ by more than $1$.
It now follows that
\[\Var_u(\E_u[Z|W]) \leq e^2\Var_u(W) < \frac{2}{5}e^2n^{-1}\]
where the second inequality follows from Lemma \ref{la}. Hence
\[\Var_u(Z) < \left(\frac{8+21e^2}{20}\right)n^{-1} < 9n^{-1}.\]
This allows us to write $Z=\E[Z|U] + D$, where $D=Z-\E[Z|U]$ has
\beqa
\Var(D) = \E\left[\Var(D|U)\right] = \E\left[\Var(Z|U)\right] < 9n^{-1}.
\eeqa
This finishes the proof of the variance part of the lemma.

Next let $n \ra \iy$, considering for each $n$ card $a_n$, where $a_n \ra a \in [0,1]$.
Then, by Lemma \ref{la},
$\E[1-W_n] \ra (1-a)e^a$ so since $\Var(W_n) \leq 2n^{-1}/5 \ra 0$, $W_n$ converges in probability to $u_0(a)=1-(1-a)e^a$.
Given $U_n=u_n \ra u$ and $W_n=w_n \ra u_0(a)$, we get by the above that for the case $u<u_0(a)$, so that $u_n<w_n$ eventually, that $Z_n \stackrel{P}{\ra} e^{1-a}u$.
For $u>u_0(a)$, so that $u_n>w_n$ eventually, $nS$ is hypergeometric and has expectation $n(u_n-w_n)(1-a_n)/(1-w_n)$.
Plugging in the limit $u_0(a)$ of $w_n$,
it follows that $S_n$ converges in probability to $1-a-e^{-a}(1-u)$.
Plugging this into (\ref{eh}) together with (\ref{ea}) and the fact that $W_n \stackrel{P}{\ra} u_0(a)$, gives that conditionally on $U=u_n \ra u$ with $u_n>W_n$,
\[Z_n \stackrel{P}{\ra} e^{e^{-a}(1-a)} - (1-u)e^{1-a}.\]
Summing up, we get that the position of a card starting from position $a_n \in Q_n$, $a_n \ra a$, after one round of CCRR shuffling converges in distribution to that of $G_a(V)$ where $V$ is uniform on $[0,1]$, as desired.
Also, taking $\mu:=\E_{u,w}[S]$, we have, since $\Var_{u,w}(S) \leq n^{-1}/4$,
\[\E_{u,w}[e^{-S}] = e^{-\mu}\E_{u,w}[e^{-(S-\mu)}] \in e^{-\mu}(1,1+\Var(S)) \subseteq e^{-\mu}\left(1 \pm \frac{1}{4n}\right).\]
Plugging this into (\ref{eh}) , letting $n \ra \iy$ and again using the convergence of $W_n$ gives (\ref{ei}).
\epf

Recall that we write $A=A(n)$ for the transition matrix of the movement of a card under one round of CCRR. Write $B=B(n)=[b_{ij}]$ for the transition matrix of a card that moves according
to $G_a(V)$. More precisely, let $V$ be uniform on $[0,1]$ and let $b_{ij}$ be the probability that $G_a(V) \in (j-1/n,j)$, $i,j \in \{1/n,2/n,\ldots,1\}$, where $a$ is chosen uniformly at random in $(i-1/n,i)$.
The precise definition of $B$ is taken so that the stationary distribution under $B$ is uniform. In particular $B$ is doubly stochastic and $B^T$ is the transition matrix of the reversed
Markov chain.

The next lemma states that the matrix $B$ has a nontrivial eigenvalue bounded away from $0$.

\bla \label{lb}
The transition matrix $B(n)$ has a (possibly complex) second eigenvalue $\lambda$ such that $|\lambda|>0.08$.
\ela

\noindent {\bf Remark.} Matlab evaluations up to $n=10^5$ strongly suggest that the second eigenvalue is real and in the interval $(0.21,0.22)$.

\smallskip

\bpf
Write $B=S+D$ where $S$ is the symmetric matrix $(B+B^T)/2$ and $D$ is the skew-symmetric matrix $(B-B^T)/2$.
We claim the following.

\bla \label{lc}
The second largest eigenvalue of $S$ is at least $0.21$
\ela

\bla \label{ld}
The (purely imaginary) eigenvalues $\lambda$ of $D$ satisfy $|\lambda|<0.13$. In particular, the
$L_{2,2}$-norm of $D$ satisfies $\n D \n_{2,2} < 0.13$.
\ela

The usefulness of Lemmas \ref{lc} and \ref{ld} and a strategy for proving them, follow from the following facts on {\em stability of eigenvalues},
i.e.\ what can happen to the spectrum of a matrix under perturbations.
These results and their elementary proofs can be found e.g.\ at \cite{TT}.
Recall that a square matrix $C$ is said to be {\em normal} if $CC^T=C^TC$ and note that $S$ and $D$ are both normal.

\bla \label{lstab}
Let $C$ be a normal $n \times n$ matrix.
Suppose that $C$ has an eigenvalue $\lambda_0$ and that $E$ is any $n \times n$ matrix with $\n E\n_{2,2} < \ep$.
Then there exists an eigenvalue $\lambda$ of $C+E$ such that $|\lambda - \lambda_0|<\ep$.

Moreover if $\lambda$ is a complex number such that there exists a vector $\phi$ such that
$\n C\phi-\lambda \phi \n_2<\ep$, then $C$ has an eigenvalue $\lambda_0$ with $|\lambda-\lambda_0|<\ep$.
\ela

Hence Lemma \ref{lb} follows immediately from Lemmas \ref{lc} and \ref{ld} together with Lemma \ref{lstab}.
\epf

\noindent {\em Proof of Lemma \ref{lc} and Lemma \ref{ld}.}
In the proof of these lemmas, it will be convenient to use the following convention: when a function $f$ is defined on $\{1/n,2/n,\ldots,1\}$ we will identify it with its extension to $[0,1]$ defined by
$f(a)=f(n^{-1}\lceil na \rceil)$.
By this convention, $\n f \n_2$ of the unextended $n$-dimensional vector $f$ is $\sqrt{n}$ times $\n f \n_2$ of the extended $f$ as a function in $L^2[0,1]$.

Let us first study $S$. That $(\lambda,\phi)$ is an eigenvalue/eigenvector pair for $S$ means that $\E[\phi(X_1)|X_0=a] = \lambda\phi(a)$ for all $a=1/n,2/n,\ldots,1$, where $X_1=X_1(n)$ is the position of a card after one move according to $S$, starting from $X_0$.
Write $Y=Y(n)$ for a random variable distributed as the position after one move according to $B(n)$ and let $Y^*(n)$ be distributed according to the position after one step of $B(n)^T$.
(Recall that $B(n)$ is doubly stochastic, so that $B(n)^T$ is the transition matrix of the reversed CCRR.)
Thus $X_1$ is the (uniform) convex combination of $Y$ and $Y^*$.
The idea now is to find $(\kappa,\psi)$ close enough to an eigenvalue/eigenvector pair to allow us to draw the desired conclusion from Lemma \ref{lstab}.
We do this with the aid of Matlab. Some more details on the Matlab computations, in particular the code, can be found in the appendix.

We use Matlab to compute the eigenvalue $\kappa=0.2293...$ and corresponding eigenvector $\chi$ with $n=10^4$, scaled so that $\n \chi \n_2=1$.
Next let $n=10^5$ and extend $\chi$ to $\psi$, the linear interpolation of (a slightly smoothed out version (see the appendix) of) $\chi$.
Then we find that
\beq \label{ej}
\n \E[\psi(X_1(n))|X_0(n)=\cdot]-\kappa\psi(\cdot))\n_2 < 0.0012.
\eeq
To arrive at the desired conclusion, a good uniform bound on the norm of the difference between $\E[\psi(X_1(m))|X_0=\cdot]$ and $\E[\psi(X_1(n))|X_0=\cdot]$ for $m > n=10^5$ will also be established.
The idea is to show that the total variation norm of the difference between the distributions of $G_a(V)$ and $G_{a+1/m}(V)$ for arbitrary $a \in (0,1)$, is small.
Then this bound will be used to infer the existence of a coupling $(Y_a,Y_{a+1/m})$ of two random variables distributed according to these, such that $\Pro(Y_a \neq Y_{a+1/m})$ small.
This together with the fact that $\n \psi \n_{\iy}$ is not too large will then establish the desired bound.

Note that the distribution function of $G_a(V)$ is $G_a^{-1}$ and the density is $(d/dx)G_a^{-1}(x)$.
Recall from (\ref{eideal}) that $u_0(a)=1-(1-a)e^a$ is the breakpoint in the expression for $G_a(u)$.

\medskip

\noindent {\em Claim.} We have
\[\n G_a(V) - G_{a+1/m}(V)\n_{TV} = \max_x|G_a^{-1}(x)-G_{a+1/m}^{-1}(x)|\]
and the difference $|G_a^{-1}(x)-G_{a+1/m}^{-1}(x)|$ is maximized when either $x=x_0:=G_{a+1/m}(u_0(a+1/m))$ or $x=G_a(u_0(a))$.

\smallskip

{\em Proof of claim.}
Write $b:=a+1/m$.
To prove the claim, it suffices to show that $(d/dx)(G_a^{-1}(x)-G_b^{-1}(x))$ is negative for $x<G_a(u_0(a))$ and $x>G_b(u_0(b))$ and positive for $G_a(u_0(a))<x<G_b(u_0(b))$.
This is equivalent to showing that $G_a'(G_a^{-1}(x))-G_b'(G_b^{-1}(x))$ is positive for $x<G_a(u_0(a))$ and $x>G_b(u_0(b))$ and negative for $x$ between the two bounds.
The derivative of $G_a$ is given by
\beq \label{eew}
G_a'(u) = \left\{ \begin{array}{ll} e^{1-a}, & u < u_0(a) \\ e^{1-a}-e^{-a}e^{e^{-a}(1-u)}, & u > u_0(a) \end{array} \right.
\eeq
Note that $G_a''\geq 0$ on $(u_0(a),1]$ so that $G_a'$ is non-decreasing on $(u_0(a),1]$.
(This is obviously true on $[0,u_0(a))$ as well, but we will not need that here.)
For $x<G_a(u_0(a))$, the difference of the derivatives is constantly $e^{1-a}-e^{1-a-1/m}>0$.
When $G_a(u_0(a)) < x < G_b(u_0(b))$,
$G_b'(G_b^{-1}(x)) = e^{1-b}$, whereas $G_a'(G_a^{-1}(x)) \leq e^{1-a}-1$, since $G_a'$ is increasing, which is obviously smaller.

For $x>G_b(u_0(b))$, let $z:=G_a^{-1}(x)$ and $y:=G_b^{-1}(x)$.
Then, since $G_a'$ is increasing, we have
\[z-y \geq \frac{G_b(z)-G_a(z)}{G_b'(z)}.\]
We want to bound this from below.
We have
\[e^{e^{-a}(1-z)}-e^{e^{-b}(1-z)} = e^{e^{-b}(1-z)}(e^{(e^{-a}-e^{-b})(1-z)} - 1).\]
Since $b > 1/m$,
\[e^{(e^{-a}-e^{-b})(1-z)} = e^{e^{-b}(e^{1/m}-1)(1-z)} \leq 1+\frac{1}{m}(1-z).\]
Hence
\beqa
\frac{G_b(z)-G_a(z)}{G_b'(z)} &=& \frac{e^{e^{-b}(1-z)}-e^{1-b}-e^{e^{-a}(1-z)}+e^{1-a}}{e^{1-b}-e^{e^{-b}(1-z)}} \\
&\geq& \frac{m^{-1}e^{1-b}-m^{-1}(1-z)e^{e^{-b}(1-z)}}{e^{1-b}-e^{e^{-b}(1-z)}} \\
&\geq& \frac{m^{-1}(1-z)(e^{1-b}-e^{e^{-b}(1-z)})}{e^{1-b}-e^{e^{-b}(1-z)}} \\
&=& m^{-1}(1-z)
\eeqa
Therefore $1-y = 1-z+z-y \geq (1+1/m)(1-z)$, so
\beqa
G_a'(z)-G_b'(y) &=& e^{1-a}-e^{1-b}+e^{e^{-b}(1-y)}-e^{e^{-a}(1-z)} \\
&\geq& m^{-1}e^{1-b} + e^{e^{-b}(1-z)(1+1/m)}-e^{e^{-a}(1-z)} \\
&\geq& m^{-1}e^{1-b} + e^{e^{-a}(1-z)} \left( e^{e^{-b}(1+1/m)-e^{-a})(1-z)}-1 \right) \\
&\geq& m^{-1}e^{1-b} + e^{e^{-a}(1-z)}\left( e^{m^{-2}e^{-a}(1-z)}-1 \right) \\
&\geq& m^{-1}e^{1-b} - 2m^{-2} \geq 0,
\eeqa
This proves the claim. \hfill $\Box$

\medskip

Now we have that
\beqa
0 &\leq& G_{a+1/m}^{-1}(G_a(u_0(a)))-G_a^{-1}(G_a(u_0(a))) \\
&=& (e^{a+1/m-1}-e^{a-1} )G_a(u_0(a)) \\
&\leq& e^{-1/m}(e^{1/m}-1)G_a(u_0(a)) \\
&\leq& m^{-1}.
\eeqa
Also
\[0 \leq G_{a}^{-1}(x_0)-G_{a+1/m}^{-1}(x_0) = G_{a}^{-1}(x_0)-u_0(a+1/m).\]
Now $G_a(u_0(a)) = e^{1-a}-e(1-a)$ from which it follows that
\beqa
x_0-G_a(u_0(a)) &=& G_{a+1/m}(u_0(a+1/m))-G_{a}(u_0(a)) \\
&=& e^{1-a}-e^{1-a-1/m}+\frac{e}{m} \\
&\leq& \frac{e}{m}\left(1-e^{-a-1/m} \right)\\
&\leq& \frac{2ea}{m}.
\eeqa
The derivative of $G_a$ was given above in (\ref{eew}) and
is minimized as $u \downarrow u_0(a)$ and then tends to $e^{1-a}-e^{1-2a} > a/2$.
Since $(d/dx)G_a^{-1}(x) = 1/G_{a}'(G_a^{-1}(x))$, it follows that
\beqa
G_a^{-1}(x_0) &\leq& u_0(a)+\frac{1}{a/2}\left(x_0-G_a(u_0(a))\right) \\
&\leq& u_0(a)+4em^{-1}.
\eeqa
Since $G_{a+1/m}^{-1}(x_0) = u_0(a+1/m)>u_0(a)$, it follows that
\[G_{a}^{-1}(x_0)-G_{a+1/m}^{-1}(x_0)<4em^{-1} < 11m^{-1}.\]
Equivalently, for $a \in \{1/n,2/n,\ldots,1-1/n\}$ and all $l$,
\begin{eqnarray}
|\sum_{j=1}^l b_{a+1/m,j}-\sum_{j=1}^l b_{a,j}| &\leq& |G_{a+1/m}^{-1}(l)-G_a^{-1}(l)| \nonumber \\
&<& 11m^{-1}. \label{eex}
\end{eqnarray}
This means that the total variation distance between the distributions of two cards making a move according to $B$, starting from $i$ and $i+1/m$ respectively, is bounded by $11/m$.
Writing $Y_a(m)$ for a random variable distributed according to the position after one round of CCRR for a card that starts in position $a$,
a consequence of this is that one can construct a coupling of $Y_{a}(m)$ and $Y_{a+1/m}(m)$ such that $\Pro(Y_a(m) \neq Y_{a+1/m}(m))<11/m$.
More generally, for $k < m$, one can couple so that $\Pro(Y_a(m) \neq Y_{a+k/m}(m))<11k/m$.
This entails, with $\hat{\psi}:=\max_x\psi(x)-\min_x\psi(x)<4.5$, that
\beq \label{el}
|\E[\psi(Y_a(m))]-\E[\psi(Y_{a+k/m}(m))]| < 11\hat{\psi}km^{-1} < 50km^{-1}.
\eeq
Next we give a bound corresponding to (\ref{el}) for $B^T$. Note that $G_a'(j)=G_a'(j+1/m)$ for $a$ such that $j+1/m<u_0(a)$ and that when $j>u_0(a)$, $G_a'(j)<G_a'(j+1/m)$, whereas when $j<u_0(a)<j+1/m$, then $G_a'(j)>G_a'(j+1/m)$.
Hence $b_{i,j+1/m}-b_{i,j}$ is zero for $u_0(i)>b+1/m$, negative for $u_0(i)<b$ and positive for the $i$'s such that $j<u_0(i)<j+1/m$.
Hence the sum
\[T:=\sum_{i:j<u_0(i)<j+1/m}(b_{i,j+1/m}-b_{i,j})\]
gives the total variation distance between the distributions of two cards making one move according to $B^T$ and starting
from $j$ and $j+1/m$ respectively.
The number of $i$'s in the sum equals at most $m(u_0^{-1}(j+1/m)-u_0^{-1}(j))+1$ and
\[b_{i,j+1/m} \leq 1 \wedge \frac{1}{m(G_i)_+'(u_0(i))} < 1 \wedge \frac{2}{mi}\]
where the second inequality follows from the bound $(G_i)_+'(u)>i/2$ from above.
Each of the $i$'s in the sum $T$ is an $i$  such that $j < u_0(i)$, i.e.\ $i > u_0^{-1}(j)$. Hence $T$ is bounded by
$2(u_0^{-1}(j+1/m)-u_0^{-1}(j)+1/m)/u_0^{-1}(j)$.
Now $u_0'(a) = ae^a$, so by the Mean Value Theorem, for some $a>u_0^{-1}(j)$,
\[u_0^{-1}(j+1/m)-u_0^{-1}(j) < \frac{e^{-a}}{ma} \leq \frac{1}{ma} < \frac{1}{u_0^{-1}(j)}.\]
Hence $T$ is bounded by $1 \wedge 2/mu_{0}^{-1}(j)^2$.
Since $u_0(a) = 1-(1-a)e^a \leq a^2$, we have $u_0^{-1}(j) \geq \sqrt{j}$
\[T \leq 1 \wedge \frac{2}{mj}.\]

Now, in analogy with the above, let $Y^*_a(m)$ be distributed as the position of a card after one move according to $B^T$, started from $a$.
Then one can construct a coupling such that $\Pro(Y^*_a(m) \neq Y^*_{a+k/m}(m)) < 1 \wedge 2k/ma$ and hence
\beq \label{em}
|\E[\psi(Y^*_a(m))]-\E[\psi(Y^*_{a+k/m}(m))]| < \hat{\psi}\left(1 \wedge \frac{2k}{ma}\right) < 4.5\left(1 \wedge \frac{2k}{ma}\right).
\eeq
Now compare $\E[\psi(Y_a(m)]$ and $\E[\psi(Y_a(n))]$.
For convenience, assume that $n|m$ and set $m=nl$. For $a=k/n-r/m$, $0 \leq r \leq l-1$, we have by convention that $\E[\psi(Y_a(n)] = \E[\psi(Y_{a_0}(n))]$, where $a_0:=n^{-1}\lceil na \rceil = (k+1)/n$.
Then (\ref{el}) shows that
\[|\E[\psi(Y_a(m))]-\E[\psi(Y_{a_0}(m))| < \frac{50}{n} = 0.0005.\]
From our Matlab calculations, we get $\max_x|\psi'(x)|<100$.
Then it is clear that
\[|\E[\psi(Y_{a_0}(m))]-\E[\psi(Y_{a_0}(n))]| \leq \frac{100}{n} = 0.001.\]
Hence
\beq \label{eo}
\n \E[\psi(Y_{\cdot}(m))]-\E[\psi(Y_{\cdot}(n))] \n_2 < 0.0015.
\eeq
Analogously for comparing $\E[\psi(Y^*_a(m)]$ with $\E[\psi(Y^*_a(n))]$, use (\ref{em}) to get
\[|\E[\psi(Y^*_a(m))]-\E[\psi(Y^*_{a_0}(m))]| < 4.5\left(1 \wedge \frac{2r}{(k-1)l}\right)\]
and hence some straightforward calculations give, using (\ref{em}), that $\sum_1^l k^2 \leq (l+1)^3/3$ and that $\sum_1^\iy 1/k^2 = \pi^2/6$,
\begin{eqnarray}
\n \E[\psi(Y^*_{\cdot}(m))]-\E[\psi(Y^*_{\cdot}(n))] \n_2 &<& 0.001+\frac{2 \cdot 4.5}{\sqrt{n}}\sqrt{1+\frac{5}{24}+\frac13\left(\frac{\pi^2}{6}-1\right)} \nonumber \\
 &<& 0.001+\frac{10}{\sqrt{n}} < 0.033 \label{eez}.
\end{eqnarray}
Since $X_1$ is the convex combination of $Y$ and $Y^*$, it follows from (\ref{eo}) and (\ref{eez}) that
\beq \label{eq}
\E[\psi(X_{\cdot}(m))]-\E[\psi(X_{\cdot}(n))] \n_2 < 0.018.
\eeq
Combining (\ref{eq}) with (\ref{ej}), we find that
\[\n \E[\psi(X_{\cdot}(m))]- \kappa\psi \n_2 < 0.0192\]
for all $m \geq 10^5$.
From this it follows that $S$ has an eigenvalue $\lambda$ with $\lambda>\kappa-0.0192>0.21$ as desired.

Next we prove Lemma \ref{ld} in a completely analogous way. We have that $(\lambda,\phi)$ is an eigenvalue/eigenvector pair for $D$ if
$(1/2)(\E[\phi_i(Y)]-\E[\phi(Y_i^*)]) = \lambda\phi(i)$ for all $i$, where $Y_i$ and $Y_i^*$ are, as above, random variables distributed according one step of $B$ and $B^T$ respectively, starting from $X_0=i$.
Again we take $n=10^5$ and use Matlab to get $\kappa$ and $\psi$ close to an eigenvalue and eigenvector respectively.
It turns out that $\kappa = 0.0793...i$, so $|\kappa|<0.08$ and we get $\hat{\psi}<5$.
In terms of variability however, this case turns out to be less well behaved. We get $\max_x|\psi'(x)|<400$ and
\[\left\n \frac12\left(\E[\psi(Y_{\cdot}(n))]-\E[\psi(Y_{\cdot}^*(n))]\right) - \kappa \psi(\cdot) \right\n_2<0.017.\]
Then the above calculations now give
\[\left\n \frac12\left(\E[\psi(Y_{\cdot}(m))]-\E[\psi(Y_{\cdot}^*(m))]\right) - \kappa \psi(\cdot)\right\n_2<0.047 < 0.05.\]
The desired result follows now follows from Lemma \ref{lstab}.
\hfill $\Box$

\medskip

For the remainder of the paper, in the light of Lemma \ref{lb}, we fix $\lambda$ to be the eigenvalue of $B$ with the second largest modulus.
Let $\phi$ be an eigenvector corresponding to $\lambda$ with $\n\phi \n_2=1$.
Note that since $\lambda$ may be complex, so may $\phi$. (However, as remarked before, Matlab computations up to $n=10^5$ strongly suggest that $\lambda$ is real.)
The next lemma, which we extract from (\ref{eex}) in the proof of Lemma \ref{lc}, will be useful in order to show that $\phi(i)$ and $\phi(j)$ cannot differ much if $i$ and $j$ are close.

\bla \label{le}
Let $f:Q_n \ra \CC$ and for $i \in Q_n$, let $X_i$ be a random variable distributed according to the law the position of card $i$ after one move according to $B$.
Then for all $i \in Q_n \setminus \{1\}$,
\[\left|\E[f(X_{i+1/n})]-\E[f(X_i)] \right| \leq \frac{22}{n}\n f\n_{\iy}.\]
\ela


%

\bla \label{lf}
For the eigenvector $\phi$, of $B$, we have
\[\n \phi\n_1 \geq c_1n^{4/9},\]
\[\n \phi\n_\iy \leq c_2n^{-4/9}\]
and
\[|\phi(i+1/n)-\phi(i)| \leq c_3n^{-13/9}\]
for constants $c_1$, $c_2$ and $c_3$ independent of $n$ and $i$.
\ela

\bpf
Let, as in Lemma \ref{le}, $X_i$ be distributed as the position of card $i$ after one move according to $B$.
By definition of eigenvalue/eigenvector, $\E[\phi(X_i)] = \lambda\phi(i)$.
Hence by Lemmas \ref{le} and \ref{lb},
\[|\phi(i)-\phi(i+1/n)| \leq \frac{22|\lambda|^{-1}}{n} < \frac{275}{n}\]
since $\n \phi \n_{\iy} \leq 1$.
Write $\n \phi\n_\iy = 100n^{-a}$.
Since $|\phi(i)-\phi(i+1/n)| < 275/n$, it follows that
\[1 \geq \n\phi\n_2^2 > \frac{275^2}{n^2}\sum_1^{\frac{100n^{1-a}}{275}} j^2  > n^{1-3a}\]
which entails that $a \geq 1/3$.
This however means that $\n\phi\n_{\iy} \leq 100n^{-1/3}$ so that the conclusion from Lemmas \ref{le} and \ref{lb} above can be strengthened to
\[|\phi(i)-\phi(i+1/n)| < \frac{27500}{n^{4/3}}.\]
Now writing $\n\phi\n_\iy = 100n^{-b}$ gives that
\[1 \geq \n\phi\n_2^2 > \frac{27500^2}{n^{8/3}}\sum_{j=1}^{\frac{100n^{1-b}}{27500}} j^2  > n^{4/3-3b}\]
so that $b \geq 4/9$.
This shows that $\n \phi\n_\iy \leq 100n^{-4/9}$.
Once again bootstrapping the bound on $|\phi(i)-\phi(i+1/n)|$ gives an upper bound of $2750000n^{-13/9}$.
Since $\n\phi\n_\iy \geq n^{-1/2}$, it follows that
\[\n\phi\n_1 = 2750000n^{-13/9}\sum_{j=1}^{\frac{n^{17/18}}{2750000}}j = n^{4/9}/5500000.\]
\epf

Let $S_t:=\sum_{i:\Re\phi(i)>0}\phi(X^i_t)$ where $X^i_t$ is the position of card $i$ after $t$ rounds of CCRR.
The random variable $S_t$ is going to be the test statistic used to verify that order $\log n$ rounds are necessary for the deck to mix.
Let $X_\iy$ be the deck at stationarity (i.e.\ uniform on $\SG_n$) and let $S_\iy=\sum_{\Re\phi(i)>0}\phi(X^i_\iy)$.
Note that $|S_0| \geq C_1n^{4/9}$ for a constant $C_1$ independent of $n$ by Lemma \ref{lf}.
Since the cards now move according to $A$ and not $B$, $\phi$ is not quite an eigenvector for the motion of a card.
However, letting $Y^i_t$ be the position of a card after $t$ steps according to $B$ and coupling $X^i_1$ and $Y^i_1$ by using the same uniform random variable $V$ for updating (for $X^i_t$ this is to say that we use $n^{-1}\lceil nV \rceil$),
(\ref{eh}) gives that $|\E[X^i_1]-\E[Y^i_1]| \leq 4/n$. Hence by Lemma \ref{lf},
\[|\E[\phi(X^i_1)|X^i_0=a]-\E[\phi(Y^i_t)|Y^i_0=a]| \leq C_2n^{-13/9}\]
for a constant $C_2$ independent of $n$.
Hence summing over $i$ with $\Re\phi(X^i_0)>0$ and using the triangle inequality gives
\[|\E[S_1|X_0]-\lambda S_0| < C_2n^{-4/9}.\]
A straightforward recursion gives, using Lemma \ref{lf},
\beq \label{en}
\E[|S_t|] \geq |\lambda|^t S_0-\left(\sum_{r=0}^{t-1}|\lambda|^r \right)C_2n^{-4/9} > C_3 |\lambda|^t n^{4/9} - C_4n^{-4/9}
\eeq
for constants $C_3$ and $C_4$ independent of $n$.

We also need to bound the variance of $S_t$.
Let $f_i(U_i) = \E[X_1^1|U_i]$, where $U_i$ is the position where card $i$ is reinserted in round 1.
Then we can write $X_1^i=f_i(U_i)+\ep_i$, where $\E[\ep_i^2] \leq C_5/n$ by (\ref{ea}).
Hence $\phi(X_i)=\phi(f_i(U_i))+\delta_i$, where the variance of $\delta_i$ is bounded by $C_6n^{-1} \cdot (n^{-4/9})^2 = C_6n^{-17/9}$ since by
Lemma \ref{lf}, $|\delta_i| \leq c_3 n^{-4/9}|\ep_i|$.

Now observe that $f_i(U_i)$ and $f_j(U_j)$ are independent and $\Cov(\delta_i,\delta_j) \leq C_6n^{-17/9}$.
Also, for all $u$ and $v$,
\begin{equation} \label{exxx}
|\E[\ep_j|U_i=u]-\E[\ep_j|U_i=v]| \leq \frac{e^2}{n}.
\end{equation}
To see this, couple the motion of card $j$ under $U_i=u$ with the motion of $j$ under $U_i=v$ by using the same $U_k$ for all $k \neq i$.
If $i<j$, then after card $j$ has been reinserted (at the same position in the two decks), the number of cards of those that remain to be reinserted that are below $j$, will differ by at most one between the two decks.
The same goes in the case $i>j$ after $i$ is reinserted.
Now use Lemma \ref{la}.

We get
\beqa
\Cov(\phi(f_i(U_i)),\delta_j) &=& \Cov(\phi(f_i(U_i)),\E[\delta_j|U_i]) \\
&\leq& \Var(\phi(f_i(U_i)))^{1/2}\Var(\E[\delta_j|U_i])^{1/2} \\
&\leq& C_7n^{-4/9}n^{-13/9} \\
&=& C_7n^{-17/9},
\eeqa
where the second inequality uses (\ref{exxx}) to bound the second factor.
Summing up, we get
\beq \label{eb}
\Cov(\phi(X_i),\phi(X_j)) \leq C_8n^{-17/9}
\eeq
from which it follows that
\beq \label{ec}
\Var(S_1) \leq C_8n^{1/9}.
\eeq
From the considerations leading up to (\ref{en}), we can write $\E[S_{t+1}|X_t] = \lambda S_t+Z$ for a random variable $Z$,  which is function of $S_t$ such that $|Z| \leq C_2n^{-4/9}$.
Hence, since $\E[|S_t|] \leq C_9n^{5/9}$ by Lemma \ref{lf}
\beqa
\Var(\E[S_{t+1}|X_t]) &\leq& |\lambda|^2\Var(S_t)+2C_2|\lambda|\E[|S_t|]n^{-4/9}+C_2^2n^{-8/9} \\ &<& |\lambda|^2\Var(S_t)+2C_2|\lambda|n^{1/9}+C_2^2n^{-8/9}.
\eeqa
Here, the mid term follows from $\Cov(S_t,Z) \leq \E[|S_t||Z|]$.
By (\ref{ec}), $\max_x\E[\Var(S_{t+1}|X_t=x)] \leq C_8n^{1/9}$.
Hence, with $v_t:=\Var(S_t)$and using that $\Var(S_{t+1}) = \E[\Var(S_{t+1}|X_t)] + \Var(\E[S_{t+1}|X_t])$, we have the recursive inequality,
\[v_{t+1} \leq |\lambda|^2 v_t+C_{10}n^{1/9}\]
with $v_0=0$.
It follows that
\[v_t \leq C_{10}n^{1/9}\sum_{j=0}^t |\lambda|^{2j} < C_{11}n^{1/9}.\]
By continuity we also get $\Var(S_\iy) \leq C_{11}n^{1/9}$.

Finally let $\tau:=\lfloor \log n/9\log |\lambda|^{-1} \rfloor$.
Then by (\ref{en}), $\E[|S_{\tau}|] \geq C_{12} n^{1/3}$, so by Chebyshev's inequality,
\[\Pro(|S_\tau| \leq n^{2/9}) \ra 0\]
as $n \ra \iy$, whereas, since $\E[S_\iy]=0$,
\[\Pro(|S_\iy| \leq n^{2/9}) \ra 1.\]
This proves the main theorem.

\section{Appendix}
For the Matlab computations, we have used three functions, {\tt rimatris}, {\tt riprod} and {\tt riprod2}.
Recall from Lemmas \ref{lb}, \ref{lc} and \ref{ld}, the transition matrix $B(n)$ for which there was established that the second eigenvalue has modulus at least $0.08$, via considerations of approximate eigenvalues and eigenvectors for the matrices $S(n)=(B(n)+B(n)^T)/2$ and $D(n)=(B(n)-B(n)^T)/2$.
The command {\tt rimatris(n)} produces $B(n)$.
The two other functions take an $n$-dimensional vector ${\mathbf v}$ as input and return $S(n){\mathbf v}$ and $D(n){\mathbf v}$ respectively.
Since we needed $n$ to be as large as $10^5$, computation time was an important issue. Therefore the code has been optimized for computational speed and it is not quite as straightforward as one would at first believe on knowing $B(n)$.
Here is the code.

\medskip
{\tt
\noindent function A=rimatris(n) \\

\smallskip

\noindent A=zeros(n,n); \\
r=zeros(1,n+1); \\
e=exp(1); \\
ep=1/n; \\
a=0; \\
ea=1; \\
ema=1; \\
eema=exp(1); \\
eep=exp(ep); \\
emep=1/eep; \\

\smallskip

\noindent for i=1:n, \\
    a=a+ep; \\
    ea=ea*eep; \\
    ema=ema*emep; \\
    eema=eema\^emep; \\
    u=0; \\

\smallskip

\noindent    for j=0:n, \\
        z=j*ep; \\
        s=min(e*ema*u,eema\^(1-u)-e*ema*(1-u))-z; \\
        while abs(s)>1e-12, \\
            I=(u <= 1-(1-a)*ea); \\
            u=u-s/(I*e*ema + (1-I)*(e*ema-ema*eema\^(1-u))); \\
            s=min(e*ema*u,eema\^(1-u)-e*ema*(1-u))-z; \\
        end \\
        r(j+1)=u; \\
    end \\
    A(i,:)=r(2:n+1)-r(1:n); \\
end
}

\medskip

{\tt
\noindent function y=riprod(x); \\

\smallskip

\noindent n=length(x); \\
y=zeros(1,n); \\
e=exp(1); \\
ep=1/n; \\
z=ep*(0:n); \\
u=z; \\
a=0; \\
ea=1; \\
ema=1; \\
eema=exp(1); \\
eep=exp(ep); \\
emep=1/eep; \\

\smallskip

\noindent for i=1:n, \\
    a=a+ep; \\
    ea=ea*eep; \\
    ema=ema*emep; \\
    eema=eema\^emep; \\

    \smallskip

\noindent  s=min(e*ema*u,eema.\^(1-u)-e*ema*(1-u))-z; \\
    while max(abs(s))>1e-12, \\
        I=(u <= 1-(1-a)*ea); \\
        u=u-s./(I*e*ema + (1-I).*(e*ema-ema*eema.\^(1-u))); \\
        s=min(e*ema*u,eema.\^(1-u)-e*ema*(1-u))-z; \\
    end \\
    r=u(2:n+1)-u(1:n); \\
    y(i)=y(i)+r*x; \\
    y=y+x(i)*r; \\
end \\
y=0.5*y';
}

\medskip

{\tt
\noindent function y=riprod2(x); \\

\smallskip

\noindent n=length(x); \\
y=zeros(1,n); \\
e=exp(1); \\
ep=1/n; \\
z=ep*(0:n); \\
u=z; \\
a=0; \\
ea=1; \\
ema=1; \\
eema=exp(1); \\
eep=exp(ep); \\
emep=1/eep; \\

\smallskip

\noindent for i=1:n, \\
    a=a+ep; \\
    ea=ea*eep; \\
    ema=ema*emep; \\
    eema=eema\^emep; \\

    \smallskip

\noindent  s=min(e*ema*u,eema.\^(1-u)-e*ema*(1-u))-z; \\
    while max(abs(s))>1e-12, \\
        I=(u <= 1-(1-a)*ea); \\
        u=u-s./(I*e*ema + (1-I).*(e*ema-ema*eema.\^(1-u))); \\
        s=min(e*ema*u,eema.\^(1-u)-e*ema*(1-u))-z; \\
    end \\
    r=u(2:n+1)-u(1:n); \\
    y(i)=y(i)+r*x; \\
    y=y-x(i)*r; \\
end \\
y=0.5*y';
}

Given these functions, they have been used with the following set of commands.

\smallskip

{\tt
\noindent A=rimatris(10001); \\
B=(A+A')/2; \\
C=(A-A')/2; \\
{[u,l]}=eigs(B,2); \\
u=u(:,2); \\
u=100*u; \\
l=l(2,2); \\
{[w,k]}=eigs(C,1);
w=100*w; \\
for i=2:25, u(26-i)=u(26)-i*(u(26)-u(25));, end \\
for i=2:75, w(76-i)=w(76)-i*(w(76)-w(75));, end \\
du=10000*(u(2:10001)-u(1:10000)); \\
dw=10000*(u(2:10001)-w(1:10000)); \\
x=0:10000; \\
xx=0:0.1:10000; \\
y=interp1(x,u,xx); \\
y=y'; \\
z=interp1(x,w,xx); \\
z=conj(z'); \\
r=riprod(y)-l*y; \\
s=riprod2(z)-k*z; \\
sqrt(r'*r/100000); \\
sqrt(s'*s/100000); \\
max(abs(du)); \\
max(abs(dw)); \\
}

Then $u$ and $w$ are first the normalized eigenvectors of $B(n)$ and $D(n)$ respectively for $n=10^4$. These are then smoothed out, whereupon $y$ and $z$ are the linear interpolations of the smoothed-out vectors. The commands $max(abs(du))$ and $max(abs(dw))$ give $\hat{\phi}$ in the respective cases.


\end{document}